\input amstex
\input amsppt.sty   
%
\catcode`\@=11

\def\East#1#2{\setboxz@h{$\m@th\ssize\;{#1}\;\;$}%
 \setbox@ne\hbox{$\m@th\ssize\;{#2}\;\;$}\setbox\tw@\hbox{$\m@th#2$}%
 \dimen@\minaw@
 \ifdim\wdz@>\dimen@ \dimen@\wdz@ \fi  \ifdim\wd@ne>\dimen@ \dimen@\wd@ne \fi
 \ifdim\wd\tw@>\z@
  \mathrel{\mathop{\hbox to\dimen@{\rightarrowfill}}\limits^{#1}_{#2}}%
 \else
  \mathrel{\mathop{\hbox to\dimen@{\rightarrowfill}}\limits^{#1}}%
 \fi}
\def\West#1#2{\setboxz@h{$\m@th\ssize\;\;{#1}\;$}%
 \setbox@ne\hbox{$\m@th\ssize\;\;{#2}\;$}\setbox\tw@\hbox{$\m@th#2$}%
 \dimen@\minaw@
 \ifdim\wdz@>\dimen@ \dimen@\wdz@ \fi \ifdim\wd@ne>\dimen@ \dimen@\wd@ne \fi
 \ifdim\wd\tw@>\z@
  \mathrel{\mathop{\hbox to\dimen@{\leftarrowfill}}\limits^{#1}_{#2}}%
 \else
  \mathrel{\mathop{\hbox to\dimen@{\leftarrowfill}}\limits^{#1}}%
 \fi}
\font\arrow@i=lams1
\font\arrow@ii=lams2
\font\arrow@iii=lams3
\font\arrow@iv=lams4
\font\arrow@v=lams5
\newbox\zer@
\newdimen\standardcgap
\standardcgap=40\p@
\newdimen\hunit
\hunit=\tw@\p@
\newdimen\standardrgap
\standardrgap=32\p@
\newdimen\vunit
\vunit=1.6\p@
\def\Cgaps#1{\RIfM@
  \standardcgap=#1\standardcgap\relax \hunit=#1\hunit\relax
 \else \nonmatherr@\Cgaps \fi}
\def\Rgaps#1{\RIfM@
  \standardrgap=#1\standardrgap\relax \vunit=#1\vunit\relax
 \else \nonmatherr@\Rgaps \fi}
\newdimen\getdim@
\def\getcgap@#1{\ifcase#1\or\getdim@\z@\else\getdim@\standardcgap\fi}
\def\getrgap@#1{\ifcase#1\getdim@\z@\else\getdim@\standardrgap\fi}
\def\cgaps#1{\RIfM@
 \cgaps@{#1}\edef\getcgap@##1{\i@=##1\relax\the\toks@}\toks@{}\else
 \nonmatherr@\cgaps\fi}
\def\rgaps#1{\RIfM@
 \rgaps@{#1}\edef\getrgap@##1{\i@=##1\relax\the\toks@}\toks@{}\else
 \nonmatherr@\rgaps\fi}
\def\Gaps@@{\gaps@@}
\def\cgaps@#1{\toks@{\ifcase\i@\or\getdim@=\z@}%
 \gaps@@\standardcgap#1;\gaps@@\gaps@@
 \edef\next@{\the\toks@\noexpand\else\noexpand\getdim@\noexpand\standardcgap
  \noexpand\fi}%
 \toks@=\expandafter{\next@}}
\def\rgaps@#1{\toks@{\ifcase\i@\getdim@=\z@}%
 \gaps@@\standardrgap#1;\gaps@@\gaps@@
 \edef\next@{\the\toks@\noexpand\else\noexpand\getdim@\noexpand\standardrgap
  \noexpand\fi}%
 \toks@=\expandafter{\next@}}
\def\gaps@@#1#2;#3{\mgaps@#1#2\mgaps@
 \edef\next@{\the\toks@\noexpand\or\noexpand\getdim@
  \noexpand#1\the\mgapstoks@@}%
 \global\toks@=\expandafter{\next@}%
 \DN@{#3}%
 \ifx\next@\Gaps@@\gdef\next@##1\gaps@@{}\else
  \gdef\next@{\gaps@@#1#3}\fi\next@}
\def\mgaps@#1{\let\mgapsnext@#1\FN@\mgaps@@}
\def\mgaps@@{\ifx\next\space@\DN@. {\FN@\mgaps@@}\else
 \DN@.{\FN@\mgaps@@@}\fi\next@.}
\def\mgaps@@@{\ifx\next\w\let\next@\mgaps@@@@\else
 \let\next@\mgaps@@@@@\fi\next@}
\newtoks\mgapstoks@@
\def\mgaps@@@@@#1\mgaps@{\getdim@\mgapsnext@\getdim@#1\getdim@
 \edef\next@{\noexpand\getdim@\the\getdim@}%
 \mgapstoks@@=\expandafter{\next@}}
\def\mgaps@@@@\w#1#2\mgaps@{\mgaps@@@@@#2\mgaps@
 \setbox\zer@\hbox{$\m@th\hskip15\p@\tsize@#1$}%
 \dimen@\wd\zer@
 \ifdim\dimen@>\getdim@ \getdim@\dimen@ \fi
 \edef\next@{\noexpand\getdim@\the\getdim@}%
 \mgapstoks@@=\expandafter{\next@}}
\def\changewidth#1#2{\setbox\zer@\hbox{$\m@th#2}%
 \hbox to\wd\zer@{\hss$\m@th#1$\hss}}
\atdef@({\FN@\ARROW@}
\def\ARROW@{\ifx\next)\let\next@\OPTIONS@\else
 \DN@{\csname\string @(\endcsname}\fi\next@}
\newif\ifoptions@
\def\OPTIONS@){\ifoptions@\let\next@\relax\else
 \DN@{\options@true\begingroup\optioncodes@}\fi\next@}
\newif\ifN@
\newif\ifE@
\newif\ifNESW@
\newif\ifH@
\newif\ifV@
\newif\ifHshort@
\expandafter\def\csname\string @(\endcsname #1,#2){%
 \ifoptions@\let\next@\endgroup\else\let\next@\relax\fi\next@
 \N@false\E@false\H@false\V@false\Hshort@false
 \ifnum#1>\z@\E@true\fi
 \ifnum#1=\z@\V@true\tX@false\tY@false\a@false\fi
 \ifnum#2>\z@\N@true\fi
 \ifnum#2=\z@\H@true\tX@false\tY@false\a@false\ifshort@\Hshort@true\fi\fi
 \NESW@false
 \ifN@\ifE@\NESW@true\fi\else\ifE@\else\NESW@true\fi\fi
 \arrow@{#1}{#2}%
 \global\options@false
 \global\scount@\z@\global\tcount@\z@\global\arrcount@\z@
 \global\s@false\global\sxdimen@\z@\global\sydimen@\z@
 \global\tX@false\global\tXdimen@i\z@\global\tXdimen@ii\z@
 \global\tY@false\global\tYdimen@i\z@\global\tYdimen@ii\z@
 \global\a@false\global\exacount@\z@
 \global\x@false\global\xdimen@\z@
 \global\X@false\global\Xdimen@\z@
 \global\y@false\global\ydimen@\z@
 \global\Y@false\global\Ydimen@\z@
 \global\p@false\global\pdimen@\z@
 \global\label@ifalse\global\label@iifalse
 \global\dl@ifalse\global\ldimen@i\z@
 \global\dl@iifalse\global\ldimen@ii\z@
 \global\short@false\global\unshort@false}
\newif\iflabel@i
\newif\iflabel@ii
\newcount\scount@
\newcount\tcount@
\newcount\arrcount@
\newif\ifs@
\newdimen\sxdimen@
\newdimen\sydimen@
\newif\iftX@
\newdimen\tXdimen@i
\newdimen\tXdimen@ii
\newif\iftY@
\newdimen\tYdimen@i
\newdimen\tYdimen@ii
\newif\ifa@
\newcount\exacount@
\newif\ifx@
\newdimen\xdimen@
\newif\ifX@
\newdimen\Xdimen@
\newif\ify@
\newdimen\ydimen@
\newif\ifY@
\newdimen\Ydimen@
\newif\ifp@
\newdimen\pdimen@
\newif\ifdl@i
\newif\ifdl@ii
\newdimen\ldimen@i
\newdimen\ldimen@ii
\newif\ifshort@
\newif\ifunshort@
\def\zero@#1{\ifnum\scount@=\z@
 \if#1e\global\scount@\m@ne\else
 \if#1t\global\scount@\tw@\else
 \if#1h\global\scount@\thr@@\else
 \if#1'\global\scount@6 \else
 \if#1`\global\scount@7 \else
 \if#1(\global\scount@8 \else
 \if#1)\global\scount@9 \else
 \if#1s\global\scount@12 \else
 \if#1H\global\scount@13 \else
 \Err@{\Invalid@@ option \string\0}\fi\fi\fi\fi\fi\fi\fi\fi\fi
 \fi}
\def\one@#1{\ifnum\tcount@=\z@
 \if#1e\global\tcount@\m@ne\else
 \if#1h\global\tcount@\tw@\else
 \if#1t\global\tcount@\thr@@\else
 \if#1'\global\tcount@4 \else
 \if#1`\global\tcount@5 \else
 \if#1(\global\tcount@10 \else
 \if#1)\global\tcount@11 \else
 \if#1s\global\tcount@12 \else
 \if#1H\global\tcount@13 \else
 \Err@{\Invalid@@ option \string\1}\fi\fi\fi\fi\fi\fi\fi\fi\fi
 \fi}
\def\a@#1{\ifnum\arrcount@=\z@
 \if#10\global\arrcount@\m@ne\else
 \if#1+\global\arrcount@\@ne\else
 \if#1-\global\arrcount@\tw@\else
 \if#1=\global\arrcount@\thr@@\else
 \Err@{\Invalid@@ option \string\a}\fi\fi\fi\fi
 \fi}
\def\ds@(#1;#2){\ifs@\else
 \global\s@true
 \sxdimen@\hunit \global\sxdimen@#1\sxdimen@\relax
 \sydimen@\vunit \global\sydimen@#2\sydimen@\relax
 \fi}
\def\dtX@(#1;#2){\iftX@\else
 \global\tX@true
 \tXdimen@i\hunit \global\tXdimen@i#1\tXdimen@i\relax
 \tXdimen@ii\vunit \global\tXdimen@ii#2\tXdimen@ii\relax
 \fi}
\def\dtY@(#1;#2){\iftY@\else
 \global\tY@true
 \tYdimen@i\hunit \global\tYdimen@i#1\tYdimen@i\relax
 \tYdimen@ii\vunit \global\tYdimen@ii#2\tYdimen@ii\relax
 \fi}
\def\da@#1{\ifa@\else\global\a@true\global\exacount@#1\relax\fi}
\def\dx@#1{\ifx@\else
 \global\x@true
 \xdimen@\hunit \global\xdimen@#1\xdimen@\relax
 \fi}
\def\dX@#1{\ifX@\else
 \global\X@true
 \Xdimen@\hunit \global\Xdimen@#1\Xdimen@\relax
 \fi}
\def\dy@#1{\ify@\else
 \global\y@true
 \ydimen@\vunit \global\ydimen@#1\ydimen@\relax
 \fi}
\def\dY@#1{\ifY@\else
 \global\Y@true
 \Ydimen@\vunit \global\Ydimen@#1\Ydimen@\relax
 \fi}
\def\p@@#1{\ifp@\else
 \global\p@true
 \pdimen@\hunit \divide\pdimen@\tw@ \global\pdimen@#1\pdimen@\relax
 \fi}
\def\L@#1{\iflabel@i\else
 \global\label@itrue \gdef\label@i{#1}%
 \fi}
\def\l@#1{\iflabel@ii\else
 \global\label@iitrue \gdef\label@ii{#1}%
 \fi}
\def\dL@#1{\ifdl@i\else
 \global\dl@itrue \ldimen@i\hunit \global\ldimen@i#1\ldimen@i\relax
 \fi}
\def\dl@#1{\ifdl@ii\else
 \global\dl@iitrue \ldimen@ii\hunit \global\ldimen@ii#1\ldimen@ii\relax
 \fi}
\def\s@{\ifunshort@\else\global\short@true\fi}
\def\uns@{\ifshort@\else\global\unshort@true\global\short@false\fi}
\def\optioncodes@{\let\0\zero@\let\1\one@\let\a\a@\let\ds\ds@\let\dtX\dtX@
 \let\dtY\dtY@\let\da\da@\let\dx\dx@\let\dX\dX@\let\dY\dY@\let\dy\dy@
 \let\p\p@@\let\L\L@\let\l\l@\let\dL\dL@\let\dl\dl@\let\s\s@\let\uns\uns@}
\def\slopes@{\\161\\152\\143\\134\\255\\126\\357\\238\\349\\45{10}\\56{11}%
 \\11{12}\\65{13}\\54{14}\\43{15}\\32{16}\\53{17}\\21{18}\\52{19}\\31{20}%
 \\41{21}\\51{22}\\61{23}}
\newcount\tan@i
\newcount\tan@ip
\newcount\tan@ii
\newcount\tan@iip
\newdimen\slope@i
\newdimen\slope@ip
\newdimen\slope@ii
\newdimen\slope@iip
\newcount\angcount@
\newcount\extracount@
\def\slope@{{\slope@i=\secondy@ \advance\slope@i-\firsty@
 \ifN@\else\multiply\slope@i\m@ne\fi
 \slope@ii=\secondx@ \advance\slope@ii-\firstx@
 \ifE@\else\multiply\slope@ii\m@ne\fi
 \ifdim\slope@ii<\z@
  \global\tan@i6 \global\tan@ii\@ne \global\angcount@23
 \else
  \dimen@\slope@i \multiply\dimen@6
  \ifdim\dimen@<\slope@ii
   \global\tan@i\@ne \global\tan@ii6 \global\angcount@\@ne
  \else
   \dimen@\slope@ii \multiply\dimen@6
   \ifdim\dimen@<\slope@i
    \global\tan@i6 \global\tan@ii\@ne \global\angcount@23
   \else
    \tan@ip\z@ \tan@iip \@ne
    \def\\##1##2##3{\global\angcount@=##3\relax
     \slope@ip\slope@i \slope@iip\slope@ii
     \multiply\slope@iip##1\relax \multiply\slope@ip##2\relax
     \ifdim\slope@iip<\slope@ip
      \global\tan@ip=##1\relax \global\tan@iip=##2\relax
     \else
      \global\tan@i=##1\relax \global\tan@ii=##2\relax
      \def\\####1####2####3{}%
     \fi}%
    \slopes@
    \slope@i=\secondy@ \advance\slope@i-\firsty@
    \ifN@\else\multiply\slope@i\m@ne\fi
    \multiply\slope@i\tan@ii \multiply\slope@i\tan@iip \multiply\slope@i\tw@
    \count@\tan@i \multiply\count@\tan@iip
    \extracount@\tan@ip \multiply\extracount@\tan@ii
    \advance\count@\extracount@
    \slope@ii=\secondx@ \advance\slope@ii-\firstx@
    \ifE@\else\multiply\slope@ii\m@ne\fi
    \multiply\slope@ii\count@
    \ifdim\slope@i<\slope@ii
     \global\tan@i=\tan@ip \global\tan@ii=\tan@iip
     \global\advance\angcount@\m@ne
    \fi
   \fi
  \fi
 \fi}%
}
\def\slope@a#1{{\def\\##1##2##3{\ifnum##3=#1\global\tan@i=##1\relax
 \global\tan@ii=##2\relax\fi}\slopes@}}
\newcount\i@
\newcount\j@
\newcount\colcount@
\newcount\Colcount@
\newcount\tcolcount@
\newdimen\rowht@
\newdimen\rowdp@
\newcount\rowcount@
\newcount\Rowcount@
\newcount\maxcolrow@
\newtoks\colwidthtoks@
\newtoks\Rowheighttoks@
\newtoks\Rowdepthtoks@
\newtoks\widthtoks@
\newtoks\Widthtoks@
\newtoks\heighttoks@
\newtoks\Heighttoks@
\newtoks\depthtoks@
\newtoks\Depthtoks@
\newif\iffirstnewCDcr@
\def\dotoks@i{%
 \global\widthtoks@=\expandafter{\the\widthtoks@\else\getdim@\z@\fi}%
 \global\heighttoks@=\expandafter{\the\heighttoks@\else\getdim@\z@\fi}%
 \global\depthtoks@=\expandafter{\the\depthtoks@\else\getdim@\z@\fi}}
\def\dotoks@ii{%
 \global\widthtoks@{\ifcase\j@}%
 \global\heighttoks@{\ifcase\j@}%
 \global\depthtoks@{\ifcase\j@}}
\def\prenewCD@#1\endnewCD{\setbox\zer@
 \vbox{%
  \def\arrow@##1##2{{}}%
  \rowcount@\m@ne \colcount@\z@ \Colcount@\z@
  \firstnewCDcr@true \toks@{}%
  \widthtoks@{\ifcase\j@}%
  \Widthtoks@{\ifcase\i@}%
  \heighttoks@{\ifcase\j@}%
  \Heighttoks@{\ifcase\i@}%
  \depthtoks@{\ifcase\j@}%
  \Depthtoks@{\ifcase\i@}%
  \Rowheighttoks@{\ifcase\i@}%
  \Rowdepthtoks@{\ifcase\i@}%
  \Let@
  \everycr{%
   \noalign{%
    \global\advance\rowcount@\@ne
    \ifnum\colcount@<\Colcount@
    \else
     \global\Colcount@=\colcount@ \global\maxcolrow@=\rowcount@
    \fi
    \global\colcount@\z@
    \iffirstnewCDcr@
     \global\firstnewCDcr@false
    \else
     \edef\next@{\the\Rowheighttoks@\noexpand\or\noexpand\getdim@\the\rowht@}%
      \global\Rowheighttoks@=\expandafter{\next@}%
     \edef\next@{\the\Rowdepthtoks@\noexpand\or\noexpand\getdim@\the\rowdp@}%
      \global\Rowdepthtoks@=\expandafter{\next@}%
     \global\rowht@\z@ \global\rowdp@\z@
     \dotoks@i
     \edef\next@{\the\Widthtoks@\noexpand\or\the\widthtoks@}%
      \global\Widthtoks@=\expandafter{\next@}%
     \edef\next@{\the\Heighttoks@\noexpand\or\the\heighttoks@}%
      \global\Heighttoks@=\expandafter{\next@}%
     \edef\next@{\the\Depthtoks@\noexpand\or\the\depthtoks@}%
      \global\Depthtoks@=\expandafter{\next@}%
     \dotoks@ii
    \fi}}%
  \tabskip\z@
  \halign{&\setbox\zer@\hbox{\vrule height10\p@ width\z@ depth\z@
   $\m@th\displaystyle{##}$}\copy\zer@
   \ifdim\ht\zer@>\rowht@ \global\rowht@\ht\zer@ \fi
   \ifdim\dp\zer@>\rowdp@ \global\rowdp@\dp\zer@ \fi
   \global\advance\colcount@\@ne
   \edef\next@{\the\widthtoks@\noexpand\or\noexpand\getdim@\the\wd\zer@}%
    \global\widthtoks@=\expandafter{\next@}%
   \edef\next@{\the\heighttoks@\noexpand\or\noexpand\getdim@\the\ht\zer@}%
    \global\heighttoks@=\expandafter{\next@}%
   \edef\next@{\the\depthtoks@\noexpand\or\noexpand\getdim@\the\dp\zer@}%
    \global\depthtoks@=\expandafter{\next@}%
   \cr#1\crcr}}%
 \Rowcount@=\rowcount@
 \global\Widthtoks@=\expandafter{\the\Widthtoks@\fi\relax}%
 \edef\Width@##1##2{\i@=##1\relax\j@=##2\relax\the\Widthtoks@}%
 \global\Heighttoks@=\expandafter{\the\Heighttoks@\fi\relax}%
 \edef\Height@##1##2{\i@=##1\relax\j@=##2\relax\the\Heighttoks@}%
 \global\Depthtoks@=\expandafter{\the\Depthtoks@\fi\relax}%
 \edef\Depth@##1##2{\i@=##1\relax\j@=##2\relax\the\Depthtoks@}%
 \edef\next@{\the\Rowheighttoks@\noexpand\fi\relax}%
 \global\Rowheighttoks@=\expandafter{\next@}%
 \edef\Rowheight@##1{\i@=##1\relax\the\Rowheighttoks@}%
 \edef\next@{\the\Rowdepthtoks@\noexpand\fi\relax}%
 \global\Rowdepthtoks@=\expandafter{\next@}%
 \edef\Rowdepth@##1{\i@=##1\relax\the\Rowdepthtoks@}%
 \colwidthtoks@{\fi}%
 \setbox\zer@\vbox{%
  \unvbox\zer@
  \count@\rowcount@
  \loop
   \unskip\unpenalty
   \setbox\zer@\lastbox
   \ifnum\count@>\maxcolrow@ \advance\count@\m@ne
   \repeat
  \hbox{%
   \unhbox\zer@
   \count@\z@
   \loop
    \unskip
    \setbox\zer@\lastbox
    \edef\next@{\noexpand\or\noexpand\getdim@\the\wd\zer@\the\colwidthtoks@}%
     \global\colwidthtoks@=\expandafter{\next@}%
    \advance\count@\@ne
    \ifnum\count@<\Colcount@
    \repeat}}%
 \edef\next@{\noexpand\ifcase\noexpand\i@\the\colwidthtoks@}%
  \global\colwidthtoks@=\expandafter{\next@}%
 \edef\Colwidth@##1{\i@=##1\relax\the\colwidthtoks@}%
 \colwidthtoks@{}\Rowheighttoks@{}\Rowdepthtoks@{}\widthtoks@{}%
 \Widthtoks@{}\heighttoks@{}\Heighttoks@{}\depthtoks@{}\Depthtoks@{}%
}
\newcount\xoff@
\newcount\yoff@
\newcount\endcount@
\newcount\rcount@
\newdimen\firstx@
\newdimen\firsty@
\newdimen\secondx@
\newdimen\secondy@
\newdimen\tocenter@
\newdimen\charht@
\newdimen\charwd@
\def\outside@{\Err@{This arrow points outside the \string\newCD}}
\newif\ifsvertex@
\newif\iftvertex@
\def\arrow@#1#2{\xoff@=#1\relax\yoff@=#2\relax
 \count@\rowcount@ \advance\count@-\yoff@
 \ifnum\count@<\@ne \outside@ \else \ifnum\count@>\Rowcount@ \outside@ \fi\fi
 \count@\colcount@ \advance\count@\xoff@
 \ifnum\count@<\@ne \outside@ \else \ifnum\count@>\Colcount@ \outside@\fi\fi
 \tcolcount@\colcount@ \advance\tcolcount@\xoff@
 \Width@\rowcount@\colcount@ \tocenter@=-\getdim@ \divide\tocenter@\tw@
 \ifdim\getdim@=\z@
  \firstx@\z@ \firsty@\mathaxis@ \svertex@true
 \else
  \svertex@false
  \ifHshort@
   \Colwidth@\colcount@
    \ifE@ \firstx@=.5\getdim@ \else \firstx@=-.5\getdim@ \fi
  \else
   \ifE@ \firstx@=\getdim@ \else \firstx@=-\getdim@ \fi
   \divide\firstx@\tw@
  \fi
  \ifE@
   \ifH@ \advance\firstx@\thr@@\p@ \else \advance\firstx@-\thr@@\p@ \fi
  \else
   \ifH@ \advance\firstx@-\thr@@\p@ \else \advance\firstx@\thr@@\p@ \fi
  \fi
  \ifN@
   \Height@\rowcount@\colcount@ \firsty@=\getdim@
   \ifV@ \advance\firsty@\thr@@\p@ \fi
  \else
   \ifV@
    \Depth@\rowcount@\colcount@ \firsty@=-\getdim@
    \advance\firsty@-\thr@@\p@
   \else
    \firsty@\z@
   \fi
  \fi
 \fi
 \ifV@
 \else
  \Colwidth@\colcount@
  \ifE@ \secondx@=\getdim@ \else \secondx@=-\getdim@ \fi
  \divide\secondx@\tw@
  \ifE@ \else \getcgap@\colcount@ \advance\secondx@-\getdim@ \fi
  \endcount@=\colcount@ \advance\endcount@\xoff@
  \count@=\colcount@
  \ifE@
   \advance\count@\@ne
   \loop
    \ifnum\count@<\endcount@
    \Colwidth@\count@ \advance\secondx@\getdim@
    \getcgap@\count@ \advance\secondx@\getdim@
    \advance\count@\@ne
    \repeat
  \else
   \advance\count@\m@ne
   \loop
    \ifnum\count@>\endcount@
    \Colwidth@\count@ \advance\secondx@-\getdim@
    \getcgap@\count@ \advance\secondx@-\getdim@
    \advance\count@\m@ne
    \repeat
  \fi
  \Colwidth@\count@ \divide\getdim@\tw@
  \ifHshort@
  \else
   \ifE@ \advance\secondx@\getdim@ \else \advance\secondx@-\getdim@ \fi
  \fi
  \ifE@ \getcgap@\count@ \advance\secondx@\getdim@ \fi
  \rcount@\rowcount@ \advance\rcount@-\yoff@
  \Width@\rcount@\count@ \divide\getdim@\tw@
  \tvertex@false
  \ifH@\ifdim\getdim@=\z@\tvertex@true\Hshort@false\fi\fi
  \ifHshort@
  \else
   \ifE@ \advance\secondx@-\getdim@ \else \advance\secondx@\getdim@ \fi
  \fi
  \iftvertex@
   \advance\secondx@.4\p@
  \else
   \ifE@ \advance\secondx@-\thr@@\p@ \else \advance\secondx@\thr@@\p@ \fi
  \fi
 \fi
 \ifH@
 \else
  \ifN@
   \Rowheight@\rowcount@ \secondy@\getdim@
  \else
   \Rowdepth@\rowcount@ \secondy@-\getdim@
   \getrgap@\rowcount@ \advance\secondy@-\getdim@
  \fi
  \endcount@=\rowcount@ \advance\endcount@-\yoff@
  \count@=\rowcount@
  \ifN@
   \advance\count@\m@ne
   \loop
    \ifnum\count@>\endcount@
    \Rowheight@\count@ \advance\secondy@\getdim@
    \Rowdepth@\count@ \advance\secondy@\getdim@
    \getrgap@\count@ \advance\secondy@\getdim@
    \advance\count@\m@ne
    \repeat
  \else
   \advance\count@\@ne
   \loop
    \ifnum\count@<\endcount@
    \Rowheight@\count@ \advance\secondy@-\getdim@
    \Rowdepth@\count@ \advance\secondy@-\getdim@
    \getrgap@\count@ \advance\secondy@-\getdim@
    \advance\count@\@ne
    \repeat
  \fi
  \tvertex@false
  \ifV@\Width@\count@\colcount@\ifdim\getdim@=\z@\tvertex@true\fi\fi
  \ifN@
   \getrgap@\count@ \advance\secondy@\getdim@
   \Rowdepth@\count@ \advance\secondy@\getdim@
   \iftvertex@
    \advance\secondy@\mathaxis@
   \else
    \Depth@\count@\tcolcount@ \advance\secondy@-\getdim@
    \advance\secondy@-\thr@@\p@
   \fi
  \else
   \Rowheight@\count@ \advance\secondy@-\getdim@
   \iftvertex@
    \advance\secondy@\mathaxis@
   \else
    \Height@\count@\tcolcount@ \advance\secondy@\getdim@
    \advance\secondy@\thr@@\p@
   \fi
  \fi
 \fi
 \ifV@\else\advance\firstx@\sxdimen@\fi
 \ifH@\else\advance\firsty@\sydimen@\fi
 \iftX@
  \advance\secondy@\tXdimen@ii
  \advance\secondx@\tXdimen@i
  \slope@
 \else
  \iftY@
   \advance\secondy@\tYdimen@ii
   \advance\secondx@\tYdimen@i
   \slope@
   \secondy@=\secondx@ \advance\secondy@-\firstx@
   \ifNESW@ \else \multiply\secondy@\m@ne \fi
   \multiply\secondy@\tan@i \divide\secondy@\tan@ii \advance\secondy@\firsty@
  \else
   \ifa@
    \slope@
    \ifNESW@ \global\advance\angcount@\exacount@ \else
      \global\advance\angcount@-\exacount@ \fi
    \ifnum\angcount@>23 \angcount@23 \fi
    \ifnum\angcount@<\@ne \angcount@\@ne \fi
    \slope@a\angcount@
    \ifY@
     \advance\secondy@\Ydimen@
    \else
     \ifX@
      \advance\secondx@\Xdimen@
      \dimen@\secondx@ \advance\dimen@-\firstx@
      \ifNESW@\else\multiply\dimen@\m@ne\fi
      \multiply\dimen@\tan@i \divide\dimen@\tan@ii
      \advance\dimen@\firsty@ \secondy@=\dimen@
     \fi
    \fi
   \else
    \ifH@\else\ifV@\else\slope@\fi\fi
   \fi
  \fi
 \fi
 \ifH@\else\ifV@\else\ifsvertex@\else
  \dimen@=6\p@ \multiply\dimen@\tan@ii
  \count@=\tan@i \advance\count@\tan@ii \divide\dimen@\count@
  \ifE@ \advance\firstx@\dimen@ \else \advance\firstx@-\dimen@ \fi
  \multiply\dimen@\tan@i \divide\dimen@\tan@ii
  \ifN@ \advance\firsty@\dimen@ \else \advance\firsty@-\dimen@ \fi
 \fi\fi\fi
 \ifp@
  \ifH@\else\ifV@\else
   \getcos@\pdimen@ \advance\firsty@\dimen@ \advance\secondy@\dimen@
   \ifNESW@ \advance\firstx@-\dimen@ii \else \advance\firstx@\dimen@ii \fi
  \fi\fi
 \fi
 \ifH@\else\ifV@\else
  \ifnum\tan@i>\tan@ii
   \charht@=10\p@ \charwd@=10\p@
   \multiply\charwd@\tan@ii \divide\charwd@\tan@i
  \else
   \charwd@=10\p@ \charht@=10\p@
   \divide\charht@\tan@ii \multiply\charht@\tan@i
  \fi
  \ifnum\tcount@=\thr@@
   \ifN@ \advance\secondy@-.3\charht@ \else\advance\secondy@.3\charht@ \fi
  \fi
  \ifnum\scount@=\tw@
   \ifE@ \advance\firstx@.3\charht@ \else \advance\firstx@-.3\charht@ \fi
  \fi
  \ifnum\tcount@=12
   \ifN@ \advance\secondy@-\charht@ \else \advance\secondy@\charht@ \fi
  \fi
  \iftY@
  \else
   \ifa@
    \ifX@
    \else
     \secondx@\secondy@ \advance\secondx@-\firsty@
     \ifNESW@\else\multiply\secondx@\m@ne\fi
     \multiply\secondx@\tan@ii \divide\secondx@\tan@i
     \advance\secondx@\firstx@
    \fi
   \fi
  \fi
 \fi\fi
 \ifH@\harrow@\else\ifV@\varrow@\else\arrow@@\fi\fi}
\newdimen\mathaxis@
\mathaxis@90\p@ \divide\mathaxis@36
\def\harrow@b{\ifE@\hskip\tocenter@\hskip\firstx@\fi}
\def\harrow@bb{\ifE@\hskip\xdimen@\else\hskip\Xdimen@\fi}
\def\harrow@e{\ifE@\else\hskip-\firstx@\hskip-\tocenter@\fi}
\def\harrow@ee{\ifE@\hskip-\Xdimen@\else\hskip-\xdimen@\fi}
\def\harrow@{\dimen@\secondx@\advance\dimen@-\firstx@
 \ifE@ \let\next@\rlap \else  \multiply\dimen@\m@ne \let\next@\llap \fi
 \next@{%
  \harrow@b
  \smash{\raise\pdimen@\hbox to\dimen@
   {\harrow@bb\arrow@ii
    \ifnum\arrcount@=\m@ne \else \ifnum\arrcount@=\thr@@ \else
     \ifE@
      \ifnum\scount@=\m@ne
      \else
       \ifcase\scount@\or\or\char118 \or\char117 \or\or\or\char119 \or
       \char120 \or\char121 \or\char122 \or\or\or\arrow@i\char125 \or
       \char117 \hskip\thr@@\p@\char117 \hskip-\thr@@\p@\fi
      \fi
     \else
      \ifnum\tcount@=\m@ne
      \else
       \ifcase\tcount@\char117 \or\or\char117 \or\char118 \or\char119 \or
       \char120\or\or\or\or\or\char121 \or\char122 \or\arrow@i\char125
       \or\char117 \hskip\thr@@\p@\char117 \hskip-\thr@@\p@\fi
      \fi
     \fi
    \fi\fi
    \dimen@\mathaxis@ \advance\dimen@.2\p@
    \dimen@ii\mathaxis@ \advance\dimen@ii-.2\p@
    \ifnum\arrcount@=\m@ne
     \let\leads@\null
    \else
     \ifcase\arrcount@
      \def\leads@{\hrule height\dimen@ depth-\dimen@ii}\or
      \def\leads@{\hrule height\dimen@ depth-\dimen@ii}\or
      \def\leads@{\hbox to10\p@{%
       \leaders\hrule height\dimen@ depth-\dimen@ii\hfil
       \hfil
      \leaders\hrule height\dimen@ depth-\dimen@ii\hskip\z@ plus2fil\relax
       \hfil
       \leaders\hrule height\dimen@ depth-\dimen@ii\hfil}}\or
     \def\leads@{\hbox{\hbox to10\p@{\dimen@\mathaxis@ \advance\dimen@1.2\p@
       \dimen@ii\dimen@ \advance\dimen@ii-.4\p@
       \leaders\hrule height\dimen@ depth-\dimen@ii\hfil}%
       \kern-10\p@
       \hbox to10\p@{\dimen@\mathaxis@ \advance\dimen@-1.2\p@
       \dimen@ii\dimen@ \advance\dimen@ii-.4\p@
       \leaders\hrule height\dimen@ depth-\dimen@ii\hfil}}}\fi
    \fi
    \cleaders\leads@\hfil
    \ifnum\arrcount@=\m@ne\else\ifnum\arrcount@=\thr@@\else
     \arrow@i
     \ifE@
      \ifnum\tcount@=\m@ne
      \else
       \ifcase\tcount@\char119 \or\or\char119 \or\char120 \or\char121 \or
       \char122 \or \or\or\or\or\char123\or\char124 \or
       \char125 \or\char119 \hskip-\thr@@\p@\char119 \hskip\thr@@\p@\fi
      \fi
     \else
      \ifcase\scount@\or\or\char120 \or\char119 \or\or\or\char121 \or\char122
      \or\char123 \or\char124 \or\or\or\char125 \or
      \char119 \hskip-\thr@@\p@\char119 \hskip\thr@@\p@\fi
     \fi
    \fi\fi
    \harrow@ee}}%
  \harrow@e}%
 \iflabel@i
  \dimen@ii\z@ \setbox\zer@\hbox{$\m@th\tsize@@\label@i$}%
  \ifnum\arrcount@=\m@ne
  \else
   \advance\dimen@ii\mathaxis@
   \advance\dimen@ii\dp\zer@ \advance\dimen@ii\tw@\p@
   \ifnum\arrcount@=\thr@@ \advance\dimen@ii\tw@\p@ \fi
  \fi
  \advance\dimen@ii\pdimen@
  \next@{\harrow@b\smash{\raise\dimen@ii\hbox to\dimen@
   {\harrow@bb\hskip\tw@\ldimen@i\hfil\box\zer@\hfil\harrow@ee}}\harrow@e}%
 \fi
 \iflabel@ii
  \ifnum\arrcount@=\m@ne
  \else
   \setbox\zer@\hbox{$\m@th\tsize@\label@ii$}%
   \dimen@ii-\ht\zer@ \advance\dimen@ii-\tw@\p@
   \ifnum\arrcount@=\thr@@ \advance\dimen@ii-\tw@\p@ \fi
   \advance\dimen@ii\mathaxis@ \advance\dimen@ii\pdimen@
   \next@{\harrow@b\smash{\raise\dimen@ii\hbox to\dimen@
    {\harrow@bb\hskip\tw@\ldimen@ii\hfil\box\zer@\hfil\harrow@ee}}\harrow@e}%
  \fi
 \fi}
\let\tsize@\tsize
\def\tsizenewCDlabels{\let\tsize@\tsize}
\def\ssizenewCDlabels{\let\tsize@\ssize}
\def\tsize@@{\ifnum\arrcount@=\m@ne\else\tsize@\fi}
\def\varrow@{\dimen@\secondy@ \advance\dimen@-\firsty@
 \ifN@ \else \multiply\dimen@\m@ne \fi
 \setbox\zer@\vbox to\dimen@
  {\ifN@ \vskip-\Ydimen@ \else \vskip\ydimen@ \fi
   \ifnum\arrcount@=\m@ne\else\ifnum\arrcount@=\thr@@\else
    \hbox{\arrow@iii
     \ifN@
      \ifnum\tcount@=\m@ne
      \else
       \ifcase\tcount@\char117 \or\or\char117 \or\char118 \or\char119 \or
       \char120 \or\or\or\or\or\char121 \or\char122 \or\char123 \or
       \vbox{\hbox{\char117 }\nointerlineskip\vskip\thr@@\p@
       \hbox{\char117 }\vskip-\thr@@\p@}\fi
      \fi
     \else
      \ifcase\scount@\or\or\char118 \or\char117 \or\or\or\char119 \or
      \char120 \or\char121 \or\char122 \or\or\or\char123 \or
      \vbox{\hbox{\char117 }\nointerlineskip\vskip\thr@@\p@
      \hbox{\char117 }\vskip-\thr@@\p@}\fi
     \fi}%
    \nointerlineskip
   \fi\fi
   \ifnum\arrcount@=\m@ne
    \let\leads@\null
   \else
    \ifcase\arrcount@\let\leads@\vrule\or\let\leads@\vrule\or
    \def\leads@{\vbox to10\p@{%
     \hrule height 1.67\p@ depth\z@ width.4\p@
     \vfil
     \hrule height 3.33\p@ depth\z@ width.4\p@
     \vfil
     \hrule height 1.67\p@ depth\z@ width.4\p@}}\or
    \def\leads@{\hbox{\vrule height\p@\hskip\tw@\p@\vrule}}\fi
   \fi
  \cleaders\leads@\vfill\nointerlineskip
   \ifnum\arrcount@=\m@ne\else\ifnum\arrcount@=\thr@@\else
    \hbox{\arrow@iv
     \ifN@
      \ifcase\scount@\or\or\char118 \or\char117 \or\or\or\char119 \or
      \char120 \or\char121 \or\char122 \or\or\or\arrow@iii\char123 \or
      \vbox{\hbox{\char117 }\nointerlineskip\vskip-\thr@@\p@
      \hbox{\char117 }\vskip\thr@@\p@}\fi
     \else
      \ifnum\tcount@=\m@ne
      \else
       \ifcase\tcount@\char117 \or\or\char117 \or\char118 \or\char119 \or
       \char120 \or\or\or\or\or\char121 \or\char122 \or\arrow@iii\char123 \or
       \vbox{\hbox{\char117 }\nointerlineskip\vskip-\thr@@\p@
       \hbox{\char117 }\vskip\thr@@\p@}\fi
      \fi
     \fi}%
   \fi\fi
   \ifN@\vskip\ydimen@\else\vskip-\Ydimen@\fi}%
 \ifN@
  \dimen@ii\firsty@
 \else
  \dimen@ii-\firsty@ \advance\dimen@ii\ht\zer@ \multiply\dimen@ii\m@ne
 \fi
 \rlap{\smash{\hskip\tocenter@ \hskip\pdimen@ \raise\dimen@ii \box\zer@}}%
 \iflabel@i
  \setbox\zer@\vbox to\dimen@{\vfil
   \hbox{$\m@th\tsize@@\label@i$}\vskip\tw@\ldimen@i\vfil}%
  \rlap{\smash{\hskip\tocenter@ \hskip\pdimen@
  \ifnum\arrcount@=\m@ne \let\next@\relax \else \let\next@\llap \fi
  \next@{\raise\dimen@ii\hbox{\ifnum\arrcount@=\m@ne \hskip-.5\wd\zer@ \fi
   \box\zer@ \ifnum\arrcount@=\m@ne \else \hskip\tw@\p@ \fi}}}}%
 \fi
 \iflabel@ii
  \ifnum\arrcount@=\m@ne
  \else
   \setbox\zer@\vbox to\dimen@{\vfil
    \hbox{$\m@th\tsize@\label@ii$}\vskip\tw@\ldimen@ii\vfil}%
   \rlap{\smash{\hskip\tocenter@ \hskip\pdimen@
   \rlap{\raise\dimen@ii\hbox{\ifnum\arrcount@=\thr@@ \hskip4.5\p@ \else
    \hskip2.5\p@ \fi\box\zer@}}}}%
  \fi
 \fi
}
\newdimen\goal@
\newdimen\shifted@
\newcount\Tcount@
\newcount\Scount@
\newbox\shaft@
\newcount\slcount@
\def\getcos@#1{%
 \ifnum\tan@i<\tan@ii
  \dimen@#1%
  \ifnum\slcount@<8 \count@9 \else \ifnum\slcount@<12 \count@8 \else
   \count@7 \fi\fi
  \multiply\dimen@\count@ \divide\dimen@10
  \dimen@ii\dimen@ \multiply\dimen@ii\tan@i \divide\dimen@ii\tan@ii
 \else
  \dimen@ii#1%
  \count@-\slcount@ \advance\count@24
  \ifnum\count@<8 \count@9 \else \ifnum\count@<12 \count@8
   \else\count@7 \fi\fi
  \multiply\dimen@ii\count@ \divide\dimen@ii10
  \dimen@\dimen@ii \multiply\dimen@\tan@ii \divide\dimen@\tan@i
 \fi}
\newdimen\adjust@
\def\Nnext@{\ifN@\let\next@\raise\else\let\next@\lower\fi}
\def\arrow@@{\slcount@\angcount@
 \ifNESW@
  \ifnum\angcount@<10
   \let\arrowfont@=\arrow@i \advance\angcount@\m@ne \multiply\angcount@13
  \else
   \ifnum\angcount@<19
    \let\arrowfont@=\arrow@ii \advance\angcount@-10 \multiply\angcount@13
   \else
    \let\arrowfont@=\arrow@iii \advance\angcount@-19 \multiply\angcount@13
  \fi\fi
  \Tcount@\angcount@
 \else
  \ifnum\angcount@<5
   \let\arrowfont@=\arrow@iii \advance\angcount@\m@ne \multiply\angcount@13
   \advance\angcount@65
  \else
   \ifnum\angcount@<14
    \let\arrowfont@=\arrow@iv \advance\angcount@-5 \multiply\angcount@13
   \else
    \ifnum\angcount@<23
     \let\arrowfont@=\arrow@v \advance\angcount@-14 \multiply\angcount@13
    \else
     \let\arrowfont@=\arrow@i \angcount@=117
  \fi\fi\fi
  \ifnum\angcount@=117 \Tcount@=115 \else\Tcount@\angcount@ \fi
 \fi
 \Scount@\Tcount@
 \ifE@
  \ifnum\tcount@=\z@ \advance\Tcount@\tw@ \else\ifnum\tcount@=13
   \advance\Tcount@\tw@ \else \advance\Tcount@\tcount@ \fi\fi
  \ifnum\scount@=\z@ \else \ifnum\scount@=13 \advance\Scount@\thr@@ \else
   \advance\Scount@\scount@ \fi\fi
 \else
  \ifcase\tcount@\advance\Tcount@\thr@@\or\or\advance\Tcount@\thr@@\or
  \advance\Tcount@\tw@\or\advance\Tcount@6 \or\advance\Tcount@7
  \or\or\or\or\or \advance\Tcount@8 \or\advance\Tcount@9 \or
  \advance\Tcount@12 \or\advance\Tcount@\thr@@\fi
  \ifcase\scount@\or\or\advance\Scount@\thr@@\or\advance\Scount@\tw@\or
  \or\or\advance\Scount@4 \or\advance\Scount@5 \or\advance\Scount@10
  \or\advance\Scount@11 \or\or\or\advance\Scount@12 \or\advance
  \Scount@\tw@\fi
 \fi
 \ifcase\arrcount@\or\or\advance\angcount@\@ne\else\fi
 \ifN@ \shifted@=\firsty@ \else\shifted@=-\firsty@ \fi
 \ifE@ \else\advance\shifted@\charht@ \fi
 \goal@=\secondy@ \advance\goal@-\firsty@
 \ifN@\else\multiply\goal@\m@ne\fi
 \setbox\shaft@\hbox{\arrowfont@\char\angcount@}%
 \ifnum\arrcount@=\thr@@
  \getcos@{1.5\p@}%
  \setbox\shaft@\hbox to\wd\shaft@{\arrowfont@
   \rlap{\hskip\dimen@ii
    \smash{\ifNESW@\let\next@\lower\else\let\next@\raise\fi
     \next@\dimen@\hbox{\arrowfont@\char\angcount@}}}%
   \rlap{\hskip-\dimen@ii
    \smash{\ifNESW@\let\next@\raise\else\let\next@\lower\fi
      \next@\dimen@\hbox{\arrowfont@\char\angcount@}}}\hfil}%
 \fi
 \rlap{\smash{\hskip\tocenter@\hskip\firstx@
  \ifnum\arrcount@=\m@ne
  \else
   \ifnum\arrcount@=\thr@@
   \else
    \ifnum\scount@=\m@ne
    \else
     \ifnum\scount@=\z@
     \else
      \setbox\zer@\hbox{\ifnum\angcount@=117 \arrow@v\else\arrowfont@\fi
       \char\Scount@}%
      \ifNESW@
       \ifnum\scount@=\tw@
        \dimen@=\shifted@ \advance\dimen@-\charht@
        \ifN@\hskip-\wd\zer@\fi
        \Nnext@
        \next@\dimen@\copy\zer@
        \ifN@\else\hskip-\wd\zer@\fi
       \else
        \Nnext@
        \ifN@\else\hskip-\wd\zer@\fi
        \next@\shifted@\copy\zer@
        \ifN@\hskip-\wd\zer@\fi
       \fi
       \ifnum\scount@=12
        \advance\shifted@\charht@ \advance\goal@-\charht@
        \ifN@ \hskip\wd\zer@ \else \hskip-\wd\zer@ \fi
       \fi
       \ifnum\scount@=13
        \getcos@{\thr@@\p@}%
        \ifN@ \hskip\dimen@ \else \hskip-\wd\zer@ \hskip-\dimen@ \fi
        \adjust@\shifted@ \advance\adjust@\dimen@ii
        \Nnext@
        \next@\adjust@\copy\zer@
        \ifN@ \hskip-\dimen@ \hskip-\wd\zer@ \else \hskip\dimen@ \fi
       \fi
      \else
       \ifN@\hskip-\wd\zer@\fi
       \ifnum\scount@=\tw@
        \ifN@ \hskip\wd\zer@ \else \hskip-\wd\zer@ \fi
        \dimen@=\shifted@ \advance\dimen@-\charht@
        \Nnext@
        \next@\dimen@\copy\zer@
        \ifN@\hskip-\wd\zer@\fi
       \else
        \Nnext@
        \next@\shifted@\copy\zer@
        \ifN@\else\hskip-\wd\zer@\fi
       \fi
       \ifnum\scount@=12
        \advance\shifted@\charht@ \advance\goal@-\charht@
        \ifN@ \hskip-\wd\zer@ \else \hskip\wd\zer@ \fi
       \fi
       \ifnum\scount@=13
        \getcos@{\thr@@\p@}%
        \ifN@ \hskip-\wd\zer@ \hskip-\dimen@ \else \hskip\dimen@ \fi
        \adjust@\shifted@ \advance\adjust@\dimen@ii
        \Nnext@
        \next@\adjust@\copy\zer@
        \ifN@ \hskip\dimen@ \else \hskip-\dimen@ \hskip-\wd\zer@ \fi
       \fi	
      \fi
  \fi\fi\fi\fi
  \ifnum\arrcount@=\m@ne
  \else
   \loop
    \ifdim\goal@>\charht@
    \ifE@\else\hskip-\charwd@\fi
    \Nnext@
    \next@\shifted@\copy\shaft@
    \ifE@\else\hskip-\charwd@\fi
    \advance\shifted@\charht@ \advance\goal@ -\charht@
    \repeat
   \ifdim\goal@>\z@
    \dimen@=\charht@ \advance\dimen@-\goal@
    \divide\dimen@\tan@i \multiply\dimen@\tan@ii
    \ifE@ \hskip-\dimen@ \else \hskip-\charwd@ \hskip\dimen@ \fi
    \adjust@=\shifted@ \advance\adjust@-\charht@ \advance\adjust@\goal@
    \Nnext@
    \next@\adjust@\copy\shaft@
    \ifE@ \else \hskip-\charwd@ \fi
   \else
    \adjust@=\shifted@ \advance\adjust@-\charht@
   \fi
  \fi
  \ifnum\arrcount@=\m@ne
  \else
   \ifnum\arrcount@=\thr@@
   \else
    \ifnum\tcount@=\m@ne
    \else
     \setbox\zer@
      \hbox{\ifnum\angcount@=117 \arrow@v\else\arrowfont@\fi\char\Tcount@}%
     \ifnum\tcount@=\thr@@
      \advance\adjust@\charht@
      \ifE@\else\ifN@\hskip-\charwd@\else\hskip-\wd\zer@\fi\fi
     \else
      \ifnum\tcount@=12
       \advance\adjust@\charht@
       \ifE@\else\ifN@\hskip-\charwd@\else\hskip-\wd\zer@\fi\fi
      \else
       \ifE@\hskip-\wd\zer@\fi
     \fi\fi
     \Nnext@
     \next@\adjust@\copy\zer@
     \ifnum\tcount@=13
      \hskip-\wd\zer@
      \getcos@{\thr@@\p@}%
      \ifE@\hskip-\dimen@ \else\hskip\dimen@ \fi
      \advance\adjust@-\dimen@ii
      \Nnext@
      \next@\adjust@\box\zer@
     \fi
  \fi\fi\fi}}%
 \iflabel@i
  \rlap{\hskip\tocenter@
  \dimen@\firstx@ \advance\dimen@\secondx@ \divide\dimen@\tw@
  \advance\dimen@\ldimen@i
  \dimen@ii\firsty@ \advance\dimen@ii\secondy@ \divide\dimen@ii\tw@
  \multiply\ldimen@i\tan@i \divide\ldimen@i\tan@ii
  \ifNESW@ \advance\dimen@ii\ldimen@i \else \advance\dimen@ii-\ldimen@i \fi
  \setbox\zer@\hbox{\ifNESW@\else\ifnum\arrcount@=\thr@@\hskip4\p@\else
   \hskip\tw@\p@\fi\fi
   $\m@th\tsize@@\label@i$\ifNESW@\ifnum\arrcount@=\thr@@\hskip4\p@\else
   \hskip\tw@\p@\fi\fi}%
  \ifnum\arrcount@=\m@ne
   \ifNESW@ \advance\dimen@.5\wd\zer@ \advance\dimen@\p@ \else
    \advance\dimen@-.5\wd\zer@ \advance\dimen@-\p@ \fi
   \advance\dimen@ii-.5\ht\zer@
  \else
   \advance\dimen@ii\dp\zer@
   \ifnum\slcount@<6 \advance\dimen@ii\tw@\p@ \fi
  \fi
  \hskip\dimen@
  \ifNESW@ \let\next@\llap \else\let\next@\rlap \fi
  \next@{\smash{\raise\dimen@ii\box\zer@}}}%
 \fi
 \iflabel@ii
  \ifnum\arrcount@=\m@ne
  \else
   \rlap{\hskip\tocenter@
   \dimen@\firstx@ \advance\dimen@\secondx@ \divide\dimen@\tw@
   \ifNESW@ \advance\dimen@\ldimen@ii \else \advance\dimen@-\ldimen@ii \fi
   \dimen@ii\firsty@ \advance\dimen@ii\secondy@ \divide\dimen@ii\tw@
   \multiply\ldimen@ii\tan@i \divide\ldimen@ii\tan@ii
   \advance\dimen@ii\ldimen@ii
   \setbox\zer@\hbox{\ifNESW@\ifnum\arrcount@=\thr@@\hskip4\p@\else
    \hskip\tw@\p@\fi\fi
    $\m@th\tsize@\label@ii$\ifNESW@\else\ifnum\arrcount@=\thr@@\hskip4\p@
    \else\hskip\tw@\p@\fi\fi}%
   \advance\dimen@ii-\ht\zer@
   \ifnum\slcount@<9 \advance\dimen@ii-\thr@@\p@ \fi
   \ifNESW@ \let\next@\rlap \else \let\next@\llap \fi
   \hskip\dimen@\next@{\smash{\raise\dimen@ii\box\zer@}}}%
  \fi
 \fi
}
\def\outnewCD@#1{\def#1{\Err@{\string#1 must not be used within \string\newCD}}}
\newskip\prenewCDskip@
\newskip\postnewCDskip@
\prenewCDskip@\z@
\postnewCDskip@\z@
\def\prenewCDspace#1{\RIfMIfI@
 \onlydmatherr@\prenewCDspace\else\advance\prenewCDskip@#1\relax\fi\else
 \onlydmatherr@\prenewCDspace\fi}
\def\postnewCDspace#1{\RIfMIfI@
 \onlydmatherr@\postnewCDspace\else\advance\postnewCDskip@#1\relax\fi\else
 \onlydmatherr@\postnewCDspace\fi}
\def\predisplayspace#1{\RIfMIfI@
 \onlydmatherr@\predisplayspace\else
 \advance\abovedisplayskip#1\relax
 \advance\abovedisplayshortskip#1\relax\fi
 \else\onlydmatherr@\prenewCDspace\fi}
\def\postdisplayspace#1{\RIfMIfI@
 \onlydmatherr@\postdisplayspace\else
 \advance\belowdisplayskip#1\relax
 \advance\belowdisplayshortskip#1\relax\fi
 \else\onlydmatherr@\postdisplayspace\fi}
\def\PrenewCDSpace#1{\global\prenewCDskip@#1\relax}
\def\PostnewCDSpace#1{\global\postnewCDskip@#1\relax}
\def\newCD#1\endnewCD{%
 \outnewCD@\cgaps\outnewCD@\rgaps\outnewCD@\Cgaps\outnewCD@\Rgaps
 \prenewCD@#1\endnewCD
 \advance\abovedisplayskip\prenewCDskip@
 \advance\abovedisplayshortskip\prenewCDskip@
 \advance\belowdisplayskip\postnewCDskip@
 \advance\belowdisplayshortskip\postnewCDskip@
 \vcenter{\vskip\prenewCDskip@ \Let@ \colcount@\@ne \rowcount@\z@
  \everycr{%
   \noalign{%
    \ifnum\rowcount@=\Rowcount@
    \else
     \global\nointerlineskip
     \getrgap@\rowcount@ \vskip\getdim@
     \global\advance\rowcount@\@ne \global\colcount@\@ne
    \fi}}%
  \tabskip\z@
  \halign{&\global\xoff@\z@ \global\yoff@\z@
   \getcgap@\colcount@ \hskip\getdim@
   \hfil\vrule height10\p@ width\z@ depth\z@
   $\m@th\displaystyle{##}$\hfil
   \global\advance\colcount@\@ne\cr
   #1\crcr}\vskip\postnewCDskip@}%
 \prenewCDskip@\z@\postnewCDskip@\z@
 \def\getcgap@##1{\ifcase##1\or\getdim@\z@\else\getdim@\standardcgap\fi}%
 \def\getrgap@##1{\ifcase##1\getdim@\z@\else\getdim@\standardrgap\fi}%
 \let\Width@\relax\let\Height@\relax\let\Depth@\relax\let\Rowheight@\relax
 \let\Rowdepth@\relax\let\Colwdith@\relax
}
\catcode`\@=\active
\hsize 30pc
\vsize 47pc
\def\nmb#1#2{#2}         
\def\cit#1#2{\ifx#1!\cite{#2}\else#2\fi} 
\def\idx{}               
\def\ign#1{}             

\redefine\o{\circ}
\define\X{\frak X}
\define\al{\alpha}
\define\be{\beta}
\define\ga{\gamma}

\define\ka{\kappa}

\define\ph{\varphi}

\define\ps{\psi}

\define\Ga{\Gamma}

\define\Om{\Omega}
\redefine\i{^{-1}}
\define\row#1#2#3{#1_{#2},\ldots,#1_{#3}}
\define\x{\times}
\define\Id{\operatorname{Id}}
\define\Fl{\operatorname{Fl}}
\define\geo{\operatorname{geo}}
\define\Vl{\operatorname{Vl}}
\define\vl{\operatorname{vl}}
\define\vpr{\operatorname{vrp}}
\def\today{\ifcase\month\or
 January\or February\or March\or April\or May\or June\or
 July\or August\or September\or October\or November\or December\fi
 \space\number\day, \number\year}
\topmatter
\title  The Jacobi Flow
\endtitle
\author  Peter W. Michor  \endauthor
\affil
Erwin Schr\"odinger Institut f\"ur Mathematische Physik,
Pasteurgasse 6/7, A-1090 Wien, Austria
\endaffil
\address
Institut f\"ur Mathematik, Universit\"at Wien,
Strudlhofgasse 4, A-1090 Wien, Austria
\endaddress
\email Peter.Michor\@esi.ac.at \endemail
\dedicatory 
For Wlodek Tulczyjew, on the occasion of his 65th birthday.
\enddedicatory
\date {\today} \enddate
\thanks 
Supported by `Fonds zur F\"orderung der wissenscahftlichen 
Forschung, Projekt P~10037~PHY'.
\endthanks
\endtopmatter

\document

It is well known that the geodesic flow on the tangent 
bundle is the flow of a certain vector field which is called the 
spray $S:TM\to TTM$. It is maybe less well known that the flow lines 
of the vector field $\ka_{TM}\o TS:TTM\to TTTM$ project to Jacobi 
fields on $TM$. This could be called the `Jacobi flow'. This result 
was developed for the lecture course \cit!{5}, and it is the main 
result of this paper.  
I was motivated by the paper \cit!{6} of Urbanski in these 
proceedings to publish it, as an explanation of some of the uses of 
iterated tangent bundles in differential geometry.

\subheading{\nmb0{1}. The tangent bundle of a vector bundle} Let
$(E,p,M)$ be a vector bundle with fiber addition 
$+_E:E\x_ME\to E$ and fiber scalar multiplication $m^E_t:E\to E$.
Then $(TE,\pi_E,E)$, the tangent bundle of the manifold $E$, 
is itself a vector bundle, with fiber addition denoted by $+_{TE}$ and
scalar multiplication denoted by $m^{TE}_t$.

If $(U_\al,\ps_\al:E\restriction U_\al\to U_\al\x V)_{\al\in A}$ 
is a vector bundle atlas for $E$, such that 
$(U_\al,u_\al)$ is a manifold atlas for $M$, then 
$(E\restriction U_\al, \ps'_\al)_{\al\in A}$ is an atlas for the 
manifold $E$, where 
$$\ps'_\al:= (u_\al\x \Id_V)\o\ps_\al:
E\restriction U_\al\to U_\al\x V \to u_\al(U_\al)\x V\subset
\Bbb R^m\x V.$$
Hence the family 
$(T(E\restriction U_\al),T\ps'_\al: T(E\restriction U_\al)\to
T(u_\al(U_\al)\x V)= u_\al(U_\al)\x V\x\Bbb R^m\x V)_{\al\in A}$ 
is the atlas describing the canonical vector bundle structure of
$(TE,\pi_E,E)$. The transition functions are in turn:
$$\align
(\ps_\al\o \ps_\be\i)(x,v)&=(x,\ps_{\al\be}(x)v) 
     \quad\text{ for }x\in U_{\al\be}\\
(u_\al\o u_\be\i)(y) &= u_{\al\be}(y) \quad\text{ for }y\in u_\be(U_{\al\be})\\
(\ps'_\al\o (\ps'_\be)\i)(y,v)&=(u_{\al\be}(y),\ps_{\al\be}(u_\be\i(y))v) \\
(T\ps'_\al\o T(\ps'_\be)\i)(y,v;\xi,w) &= 
	 \bigl(u_{\al\be}(y),\ps_{\al\be}(u_\be\i(y))v;d(u_{\al\be})(y)\xi,\\ 
&\qquad(d(\ps_{\al\be}\o u_\be\i)(y))\xi)v + \ps_{\al\be}(u_\be\i(y))w\bigr).
\endalign$$
So we see that for fixed $(y,v)$ the transition functions are linear
in $(\xi,w)\in \Bbb R^m\x V$. This describes the vector bundle
structure of the tangent bundle $(TE,\pi_E,E)$.

For fixed $(y,\xi)$ the transition functions of $TE$ are also
linear in $(v,w)\in V\x V$. This gives a vector bundle
structure on $(TE,Tp,TM)$. Its fiber addition will be denoted by
$T(+_E): T(E\x_ME)=TE\x_{TM}TE \to TE$, since it is the tangent
mapping of $+_E$. Likewise its scalar multiplication will be
denoted by $T(m^E_t)$. One may say that the second vector bundle
structure on  $TE$, that one over $TM$, is the derivative of the
original one on $E$.

The space $\{\Xi\in TE:Tp.\Xi=0\text{ in }TM\} = (Tp)\i(0)$ is
denoted by $VE$ and is called the \idx{\it vertical bundle} over $E$.
The local form of a vertical vector $\Xi$ is 
$T\ps'_\al.\Xi = (y,v;0,w)$, so the transition
functions are   
$(T\ps'_\al\o T(\ps'_\be)\i)(y,v;0,w) =
(u_{\al\be}(y),\ps_{\al\be}(u_\be\i(y))v;0,\ps_{\al\be}(u_\be\i(y))w)$. 
They are linear in $(v,w)\in V\x V$ for fixed $y$, so $VE$ is a
vector bundle over $M$. It coincides with $0_M^*(TE,Tp,TM)$, the
pullback of the bundle $TE\to TM$ over the zero section.
We have a canonical isomorphism $\Vl_E:E\x_ME\to VE$, called the
\idx{\it big vertical lift}, given by 
$\Vl_E(u_x,v_x):={\partial_t}|_0(u_x+tv_x)$, 
which is fiber linear over $M$. 
We will mainly use the \idx{\it small vertical lift} $\vl_E:E\to TE$,
given by $\vl_E(v_x)={\partial_t}|_0 t.v_x = \Vl_E(0_x,v_x)$.
The local
representation of the vertical lift is 
$(T\ps'_\al\o \vl_E\o(\ps'_\al)\i)(y,v) = (y,0;0,v)$.

If $\ph:(E,p,M)\to (F,q,N)$ is a vector bundle
homomorphism, then we have 
$\vl_F\o \ph = T\ph\o \vl_E: E\to VF\subset TF$. 
So $\vl$ is a natural transformation between certain functors
on the category of vector bundles and their homomorphisms.
The mapping $\vpr_E:= pr_2\o \Vl_E\i:VE\to E$ is called the
\idx{\it vertical projection}. 

\subheading{\nmb0{2}. The second tangent bundle of a manifold}
All of \nmb!{1} is valid for the second tangent bundle
$TTM$ of a manifold, but here we have one more natural
structure at our disposal. The \idx{\it canonical flip} or \idx{\it
involution} $\ka_M:TTM\to TTM$ is defined locally by 
$$(TTu\o\ka_M\o TTu\i)(x,\xi;\eta,\zeta) = (x,\eta;\xi,\zeta),$$
where $(U,u)$ is a chart on $M$. Clearly this definition is
invariant under changes of charts ($Tu_\al$ equals 
$\ps'_\al$ from \nmb!{1}).

The flip $\ka_M$ has the following properties:
\roster
\item $\ka_N\o TTf=TTf\o\ka_M$ for each $f\in C^\infty(M,N)$.
\item $T(\pi_M)\o\ka_M = \pi_{TM}$ and $\pi_{TM}\o\ka_M=T(\pi_M)$.
\item $\ka_M\i=\ka_M$.
\item $\ka_M$ is a linear isomorphism from the vector bundle
        $(TTM,T(\pi_M),TM)$ to the bundle $(TTM,\pi_{TM},TM)$, so it
        interchanges the two vector bundle structures on $TTM$.
\item It is the unique smooth mapping $TTM\to TTM$ which satisfies 
$$
{\partial_t}{\partial_s}c(t,s) = \ka_M{\partial_s}{\partial_t}c(t,s)
$$
        for each $c:\Bbb R^2\to M$.
\endroster
All this follows from the local formula given above.
A quite early use of $\ka_M$ is in \cit!{4}. 

\proclaim{\nmb0{3}. Lemma} For vector fields $X$, $Y\in \X(M)$
we have 
$$\align
[X,Y] &= \vpr_{TM} \o (TY\o X - \ka_M\o TX\o Y), \\
TY\o X  -_{TM} \ka_T\o TX\o Y &= \Vl_{TM}(Y,[X,Y]) \\
&= (\vl_{TM}\o [X,Y])\;T(+_{TM})\; (0_{TM}\o Y).
\endalign$$
\endproclaim
See \cit!{3}~6.13,~6.19, or 37.13 for different proofs of this well 
known result.

\subheading{\nmb0{4}. Linear connections and their curvatures} 
Let $(E,p,M)$ be a vector bundle. Recall 
that a linear connection on the vector bundle $E$ can be 
described by specifying its \idx{\it connector} $K:TE\to E$. This 
notions seems to be due to \cit!{2}.
Any smooth mapping $K:TE\to E$ which 
is a (fiber linear) homomorphism for both vector bundle structures on $TE$,
$$\cgaps{;}\newCD
TE @()\L{\pi_E}@(0,-1) @()\L{K}@(1,0) & E @()\L{p}@(0,-1) & 
 & TE @()\L{Tp}@(0,-1) @()\L{K}@(1,0) & E @()\L{p}@(0,-1) \\
E @()\L{p}@(1,0) & M & & TM @()\L{\pi_M}@(1,0) & M 
\endnewCD$$
and which is a left inverse to the vertical lift, 
$K\o \vl_E = \Id_E:E\to TE\to E$,
specifies a linear connection. Namely: The inverse image $H:=K\i(0_E)$ 
of the zero section $0_E\subset E$, it is a subvector bundle for both 
vector bundle structures, and for the vector bundle stucture 
$\pi_E:TE\to E$ the subbundle $H$ turns out to be a complementary 
bundle for the vertical bundle $VE\to E$. We get then the associated 
\idx{\it horizontal lift mapping}
$$\gather
C: TM\x_M E \to TE, \quad C(\quad,u) = 
     \Bigl(Tp|\ker(K:T_uE\to E_{p(u)}) \Bigr)\i
\endgather$$
which has the following properties
$$\align
&(Tp,\pi_E)\o C = \Id_{TM\x_ME},\\
&C(\quad,u):T_{p(u)}M\to T_uE \text{ is linear for each } u \in E,\\
&C(X_x,\quad):E_x\to (Tp)\i(X_x) \text{ is linear for each } X_x \in T_xM.\\
\endalign$$ 
Conversely given a smooth horizontal lift mapping $C$ with these 
properties one can reconstruct a connector $K$.

For any manifold $N$, smooth mapping $s:N\to E$ along $f=p\o s:N\to M$, 
and vector field $X\in \X(N)$ a connector $K:TE\to E$
defines the 
\idx{\it covariant derivative} of $s$ along $X$ by
$$\nabla_Xs:= K\o Ts\o X: N\to TN\to TE\to E.\tag1$$
See the following diagram for all the mappings.
$$\cgaps{;}\newCD
 & TE @()\L{K}@(1,-1) @()\L{\pi_E}@(0,-1) & \\
TN @()\L{Ts}@(1,1) & E & E \\
N @()\L{X}@(0,1) @()\L{s}@(1,1) 
     @()\l{\nabla_Xs}@(2,1) @()\l{f}@(1,0) & M & 
\endnewCD\tag2$$
In canonical coordinates as in \nmb!{1} we have then
$$\align
&C((y,\xi),(y,v)) = (y,v;\xi,\Ga_y(v,\xi)),\\
&K(y,v;\xi,w) = (y,w-\Ga_y(v,\xi)),\\
&\nabla_{(y,\xi)}(\Id,s) = (\Id, ds(y)\xi-\Ga_y(s(y),\xi)),
\endalign$$
where the \idx{\it Christoffel symbol} $\Ga_y(v,\xi)$ is smooth in 
$y$ and bilinear in $(v,\xi)$. Here the sign is the negative of the 
one in many more traditional approaches, since $\Ga$ parametrizes the 
horizontal bundle. 

Let $C^\infty_f(N,E)$ denote the space of all sections along $f$ of 
$E$, isomorphic to the space $C^\infty(f^*E)$ of sections of the 
pullback bundle. 
The covariant derivative may then be viewed as a bilinear mapping 
$\nabla: \X(N)\x C^\infty_f(N,E)\to C^\infty_f(N,E)$. It has the 
following properties which follow directly from the definitions:
\roster
\item[3] $\nabla_Xs$ is $C^\infty(N,\Bbb R)$-linear in $X\in\X(N)$.
    For $x\in N$ also we have 
    $\nabla_{X(x)}s=K.Ts.X(x)=(\nabla_Xs)(x)\in E$.  
\item $\nabla_X(h.s)=dh(X).s+h.\nabla_Xs$ for 
    $h\in C^\infty(N,\Bbb R)$.
\item For any manifold $Q$, smooth mapping $g:Q\to N$, and
    $Y_y\in T_yQ$ we have  
    $\nabla_{Tg.Y_y}s = \nabla_{Y_y}(s\o g)$.
    If $Y\in \X(Q)$ and $X\in \X(N)$ are $g$-related, then 
    we have $\nabla_Y(s\o g) = (\nabla_Xs)\o g$.
\endroster

For vector fields $X$, $Y\in\X(M)$ and a section $s\in C^\infty(E)$ 
the curvature $R\in \Om^2(M,L(E,E))$ of the 
connection is given by  
$$R(X,Y)s =([\nabla_X,\nabla_Y]-\nabla_{[X,Y]})s\tag6$$

\proclaim{Theorem} 
Let $K:TE\to E$ be the connector 
of a linear connection on a vector bundle $(E,p,M)$. 
If $s:N\to E$ is a section along $f:=p\o s:N\to M$ then we have
     for vector fields $X$, $Y\in\X(N)$
$$\align
\nabla_X\nabla_Ys&-\nabla_Y\nabla_Xs-\nabla_{[X,Y]}s =\tag7\\
&=(K\o TK \o \ka_E - K\o TK)\o TTs \o TX \o Y=\\
&=R\o(Tf\o X,Tf\o Y)s:N\to E,
\endalign$$
where $R\in\Om^2(M;L(E,E))$ is the curvature.
\endproclaim

\demo{Proof}
Let first $m^E_t:E\to E$ denote the scalar 
multiplication. Then we  
have ${\partial_t}|_0 m^E_t = \vl_E$ where $\vl_E:E\to TE$ is 
the vertical lift. We use then lemma \nmb!{3} and some obvious 
commutation relations to get in turn:
$$\align
\vl_E\o K &= {\partial_t}|_0 m^E_t\o K 
     = {\partial_t}|_0 K\o m^{TE}_t 
     = TK \o {\partial_t}|_0 m^{TE}_t  
     = TK \o \vl_{(TE,\pi_E,E)}.\\
\nabla_X\nabla_Ys&-\nabla_Y\nabla_Xs-\nabla_{[X,Y]}s\\
&= K\o T(K\o Ts \o Y)\o X - K\o T(K\o Ts \o X)\o Y - K\o Ts\o [X,Y]\\
K\o T&s\o [X,Y] = K\o \vl_E\o K \o Ts \o [X,Y] \\
&= K\o TK \o \vl_{TE} \o Ts \o [X,Y] 
     = K\o TK \o TTs \o \vl_{TN} \o [X,Y] \\
&= K\o TK \o TTs \o ((TY\o X - \ka_N \o TX\o Y)\;(T-)\;0_{TN}\o Y) \\
&= K\o TK \o TTs \o TY\o X - K\o TK \o TTs \o \ka_N \o  TX\o Y - 0.
\endalign$$
Now we sum up and use $TTs\o \ka_N = \ka_E\o TTs$ to get the first
result. 
If in particular we choose $f=\Id_M$ so that $s$ is a section 
of $E\to M$ and $X,Y$ are vector fields on $M$, 
then we get the curvature $R$. 

To see that in the general case 
$(K\o TK \o \ka_E - K\o TK)\o TTs \o TX \o Y$ coincides with 
$R(Tf\o X,Tf\o Y)s$ one has to 
write out \therosteritem1 and 
$(TTs \o TX \o Y)(x)\in TTE$ in canonical charts induced from 
vector bundle charts of $E$. 
\qed\enddemo

\proclaim{\nmb0{5}. Torsion}
Let $K:TTM\to M$ be a linear connector on the tangent bundle, let 
$X,Y\in\X(M)$. 
Then the torsion is given by  
$$
\operatorname{Tor}(X,Y) = (K \o \ka_M - K) \o TX \o Y.$$
If moreover $f:N\to M$ is smooth and $U,V\in \X(N)$ then we get also
$$\align
\operatorname{Tor}(Tf.U,Tf.V)
     &=\nabla_U(Tf\o V) - \nabla_V(Tf\o U) - Tf\o [U,V] \\
&= (K \o \ka_M - K) \o TTf\o TU \o V.
\endalign$$
\endproclaim

\demo{Proof}
\therosteritem9 We have in turn
$$\align
\operatorname{Tor}(X,Y) &= \nabla_XY - \nabla_YX - [X,Y]\\
&= K \o TY \o X - K \o TX \o Y - K \o \vl_{TM}\o [X,Y] \\
K \o \vl_{TM}\o [X,Y] 
     &= K\o ((TY\o X - \ka_{M}\o TX\o Y)\;(T-)\;0_{TM}\o Y)\\
&= K \o TY \o X - K \o \ka_M \o TX \o Y - 0. 
\endalign$$
An analogous computation works in the second case, and that 
$(K \o \ka_M - K) \o TTf\o TU \o V = \operatorname{Tor}(Tf.U, Tf.V)$ 
can again be checked in local coordinates.
\qed\enddemo

\subhead\nmb0{6}. Sprays \endsubhead
Given a linear connector 
$K:TTM\to M$ on the tangent bundle with its horizontal lift mapping
$C:TM\x_MTM \to TTM$, then 
$S:= C\o \operatorname{diag}:TM \to TM\x_MTM\to TTM$ is called the 
\idx{\it spray}. This notion is due to \cit!{1}. 
The spray has the following properties:
$$\alignat3
& \pi_{TM}\o S=\Id_{TM} &\quad& 
     \text{ a vector field on } TM,\\
& T(\pi_M)\o S=\Id_{TM} &\quad& 
     \text{ a second order differential equation},\\
& S\o m^{TM}_t = T(m^{TM}_t)\o m^{TTM}_t \o S &\quad& 
     \text{ `quadratic'},
\endalignat$$
where $m_t^{E}$ is the scalar multiplication by $t$ on a vector 
bundle $E$. From $S$ one can reconstruct the torsion free part of $C$.
The following result is well known:

\proclaim{Lemma} For a spray $S:TM\to TTM$ on $M$, for $X\in TM$
$$
\geo^S(X)(t) := \pi_M(\Fl_t^{S}(X))
$$
defines a \idx{\it geodesic structure} on $M$, where $\Fl^S$ is the 
flow of the vector field $S$. 
\endproclaim
The abstract properties of a geodesic structure are obvious:
$$\align
&\geo: TM\x \Bbb R\supset U \to M\\
&\geo(X)(0) = \pi_M(X), \quad {\partial_t}|_0 \geo(X)(t)=X\\
&\geo(tX)(s) = \geo(X)(ts)\\
&\geo(\geo(X)'(t))(s) = \geo(X)(t+s)
\endalign$$
From a geodesic structure one can reconstruct the spray by 
differentiation.

\proclaim{\nmb0{7}. Theorem}
Let $S:TM\to TTM$ be a spray on a manifold $M$. Then 
$\ka_{TM}\o TS:TTM\to TTTM$ is a vector field. Consider a flow line 
$$
Y(t)=\Fl^{\ka_{TM}\o TS}_t(Y(0))
$$
of this field. Then we have:
\roster
\item"" $c:=\pi_M\o\pi_{TM}\o Y$ is a geodesic on $M$.
\item"" $\dot c=\pi_{TM}\o Y$ is the velocity field of $c$.
\item"" $J:=T(\pi_M)\o Y$ is a Jacobi field along $c$.
\item"" $\dot J = \ka_M\o Y$ is the velocity field of $J$.
\item"" $\nabla_{\partial_t}J=K\o \ka_M\o Y$ 
        is the covariant derivative of $J$. 
\item"" The Jacobi equation is given by:
$$\align
0 &= \nabla_{\partial_t}\nabla_{\partial_t}J + R(J,\dot c)\dot c +
     \nabla_{\partial_t}\operatorname{Tor}(J,\dot c)\\
&= K\o TK\o TS\o Y.
\endalign$$
\endroster
This implies that in a canonical chart induced from a chart on $M$ 
the curve $Y(t)$ is given by 
$$
(c(t),c'(t);J(t),J'(t)).
$$
\endproclaim

\demo{Proof}
Consider a curve $s\mapsto X(s)$ in $TM$. Then each 
$t\mapsto \pi_M(\Fl^S_t(X(s)))$ is a geodesic in $M$, and in the 
variable $s$ it is a variation through geodesics. Thus 
$J(t):={\partial_s}|_0 \pi_M(\Fl^S_t(X(s)))$
is a Jacobi field along the geodesic $c(t):=\pi_M(\Fl^S_t(X(0)))$, 
and each Jacobi field is of this form, for a suitable curve $X(s)$.
We consider now the curve 
$Y(t):= {\partial_s}|_0 \Fl^S_t(X(s))$
in $TTM$. Then by \nmb!{2}.(6) we have
$$\align
{\partial_t} Y(t) &= {\partial_t}{\partial_s}|_0 \Fl^S_t(X(s))
     = \ka_{TM}{\partial_s}|_0{\partial_t} \Fl^S_t(X(s))
     = \ka_{TM}{\partial_s}|_0 S(\Fl^S_t(X(s)))\\
&= (\ka_{TM}\o TS)({\partial_s}|_0\Fl^S_t(X(s)))
     = (\ka_{TM}\o TS)(Y(t)),
\endalign$$
so that $Y(t)$ is a flow line of the vector field 
$\ka_{TM}\o TS: TTM\to TTTM$. Moreover using the properties of $\ka$ 
from section \nmb!{2} and of $S$ from section \nmb!{6} we get
$$\align
T\pi_M.Y(t) &= T\pi_M.{\partial_s}|_0 \Fl^S_t(X(s))
     = {\partial_s}|_0 \pi_M(\Fl^S_t(X(s))) = J(t),\\
\pi_M T\pi_M Y(t) &= c(t),\text{ the geodesic}, \\
{\partial_t} J(t) 
     &= {\partial_t}
     T\pi_M.{\partial_s}|_0 \Fl^S_t(X(s))
     = {\partial_t}
     {\partial_s}|_0 \pi_M(\Fl^S_t(X(s))),\\
&= \ka_M{\partial_s}|_0 
     {\partial_t}\pi_M(\Fl^S_t(X(s)))
     = \ka_M{\partial_s}|_0 
     {\partial_t}\pi_M(\Fl^S_t(X(s)))\\
&= \ka_M{\partial_s}|_0 
     T\pi_M.{\partial_t}\Fl^S_t(X(s))
     = \ka_M{\partial_s}|_0 
     (T\pi_M\o S)\Fl^S_t(X(s))\\
&= \ka_M{\partial_s}|_0 \Fl^S_t(X(s))
     = \ka_M Y(t),\\
\nabla_{{\partial_t}}J 
     &= K\o {\partial_t} J = K\o \ka_M\o Y.
\endalign$$
Finally let us express the well known Jacobi expression, where we put 
$\ga(t,s):= \pi_M(\Fl^S_t(X(s)))$ for short and use most of the 
expressions from above:
$$\align
\nabla_{\partial_t}&\nabla_{\partial_t}J + R(J,\dot c)\dot c +
     \nabla_{\partial_t}\operatorname{Tor}(J,\dot c)=\\
&= \nabla_{\partial_t}\nabla_{\partial_t}.T\ga.\partial_s 
     + R(T\ga.\partial_s,T\ga.\partial_t)T\ga.\partial_t 
     + \nabla_{\partial_t}\operatorname{Tor}
     (T\ga.\partial_s,T\ga.\partial_t)\\
&= K.T(K.T(T\ga.\partial_s).\partial_t).\partial_t\\
&\quad + (K.TK.\ka_{TM}-K.TK).TT(T\ga.\partial_t).T\partial_s.\partial_t\\
&\quad + K.T((K.\ka_M - K).TT\ga.T\partial_s.\partial_t).\partial_t
\endalign$$
Note that for example for the term in the second summand we have 
$$
TTT\ga.TT\partial_t.T\partial_s.\partial_t 
= T(T(\partial_t\ga).\partial_s).\partial_t
=\partial_t\partial_s\partial_t\ga 
= \partial_t.\ka_M.\partial_t.\partial_s\ga 
= T\ka_M.\partial_t.\partial_t.\partial_s\ga 
$$
which at $s=0$ equals $T\ka_M\ddot J$. 
Using this we get for the Jacobi expression at $s=0$:
$$\align
\nabla_{\partial_t}&\nabla_{\partial_t}J + R(J,\dot c)\dot c +
     \nabla_{\partial_t}\operatorname{Tor}(J,\dot c)=\\
&=(K.TK +K.TK.\ka_{TM}.T\ka_M -K.TK.T\ka_M +K.TK.T\ka_M 
-K.TK).\partial_t\partial_t J=\\
&= K.TK.\ka_{TM}.T\ka_M.\partial_t\partial_t  J 
= K.TK.\ka_{TM}.\partial_t Y = K.TK.TS.Y,\\
\endalign$$
where we used $\partial_t\partial_t J = \partial_t (\ka_M.Y) = T\ka_M 
\partial_t Y = T\ka_M.\ka_{TM}.TS.Y$. 
Finally the validity of the Jacobi equation $0=K.TK.TS.Y$ follows 
trivially from $K\o S= 0_{TM}$. 
\qed
\enddemo


\Refs


\ref
\no \cit0{1}
\by Ambrose, W; Palais, R.S., Singer I.M.
\paper Sprays
\jour Acad. Brasileira de Ciencias
\vol 32
\yr 1960
\pages 163--178
\endref

\ref
\no \cit0{2}
\by Gromoll, D.; Klingenberg, W.; Mayer, W.
\book Riemannsche Geometrie im Gro\ss en
\bookinfo Lecture Notes in Math. 55 
\publ Springer-Verlag
\yr 1968
\publaddr Berlin, Heidelberg
\endref

\ref 
\no \cit0{3}
\by Kol\'a\v r, Ivan; Michor, Peter W.; Slov\'ak, Jan
\book Natural operations in differential geometry  
\publ Springer-Verlag
\publaddr Berlin, Heidelberg, New~York
\yr 1993
\endref

\ref
\no \cit0{4}
\by Losik, M.V.
\paper On infinitesimal connections in tangential stratifiable spaces
\lang Russian
\jour Izv. Vyssh. Uchebn. Zaved., Mat.
\vol 5(42)
\yr 1964
\pages 54--60
\endref

\ref
\no \cit0{5}
\by Michor, Peter W.
\book Riemannsche Differentialgeometrie
\bookinfo Lecture course at the Universit\"at Wien
\yr 1988/89
\endref

\ref
\no \cit0{6}
\by Urbanski, P.
\paper Double bundles
\jour These proceedings
\vol ??
\yr ??
\pages ??
\endref

\endRefs
\enddocument